\documentclass[12pt]{amsart}
\usepackage{amsfonts,latexsym,amsthm,amssymb,graphicx}
\usepackage[all]{xy}
\DeclareFontFamily{OT1}{rsfs}{}
\DeclareFontShape{OT1}{rsfs}{n}{it}{<-> rsfs10}{}
\DeclareMathAlphabet{\mathscr}{OT1}{rsfs}{n}{it}

\newtheorem{theorem}{Theorem}[section]
\newtheorem{lemma}[theorem]{Lemma}

\newtheorem{claim}[theorem]{Claim}

{\theoremstyle{definition} }
{\theoremstyle{remark} \newtheorem{remark}[theorem]{Remark}
}

\newcommand{\Cbb}{{\mathbb{C}}}
\newcommand{\Nbb}{{\mathbb{N}}}
\newcommand{\Pbb}{{\mathbb{P}}}
\newcommand{\PP}{{Q}}
\newcommand{\Sym}{\text{\rm Sym}}
\newcommand{\cI}{{\mathscr I}}
\newcommand{\cK}{{\mathscr K}}
\newcommand{\cN}{{\mathscr N}}
\newcommand{\cO}{{\mathscr O}}
\newcommand{\Til}[1]{{\widetilde{#1}}}
\newcommand{\one}{1\hskip-3.5pt1}
\newcommand{\csm}{{c_{\text{SM}}}}
\DeclareMathOperator{\Spec}{Spec}
\DeclareMathOperator{\Proj}{Proj}

\title{
Chern classes of blow-ups
}
\author{Paolo Aluffi}
\address{
Mathematics Department, 
Florida State University,
Tallahassee FL 32306, U.S.A.
}
\email{aluffi@math.fsu.edu}

\begin{document}

\begin{abstract}
We extend the classical formula of Porteous for blowing-up Chern classes 
to the case of blow-ups of possibly singular varieties along regularly 
embedded centers.
The proof of this generalization is perhaps conceptually simpler than the 
standard argument for the nonsingular case, involving Riemann-Roch 
without denominators. The new approach
relies on the explicit computation of an ideal, and a mild generalization 
of the well-known formula for the normal bundle of a proper transform 
(\cite{MR85k:14004}, B.6.10).

We also discuss alternative, very short proofs of the standard formula 
in some cases: an approach relying on the theory of
Chern-Schwartz-MacPherson classes (working in characteristic~$0$), 
and an argument reducing the formula to a straightforward computation 
of Chern classes for sheaves of differential $1$-forms with logarithmic poles
(when the center of the blow-up is a complete intersection).
\end{abstract}

\maketitle


\section{Introduction}\label{intro}
\subsection{}\label{introlead}
A general formula for the Chern classes of the tangent bundle of
the blow-up of a nonsingular variety along a nonsingular center was 
conjectured by J.~A.~Todd and B.~Segre, who established several
particular cases (\cite{MR0004487}, \cite{MR0065969}). The formula
was eventually proved by I.~R.~Porteous (\cite{MR0121813}), using 
Riemann-Roch.
F.~Hirzebruch's summary of Porteous' argument in his review
of the paper (MR0121813) may be recommend for a sharp and
lucid account. For a thorough treatment, detailing the use of 
Riemann-Roch `without denominators', the standard reference is
\S15.4 in \cite{MR85k:14004}. Here is the formula in the
notation of the latter reference. For any non-singular variety $X$,
write $c(X)$ for $c(T_X)\cap [X]$, the total Chern class (in the Chow group
of~$X$) of the tangent bundle of $X$. Let $X\subseteq Y$ be nonsingular
varieties, and let $\Til Y$ be the blow-up of $Y$ along $X$, with
exceptional divisor $\Til X$:
\[
\xymatrix@R=20pt@C=20pt{
\Til X \ar[r]^j \ar[d]_g & \Til Y \ar[d]^f \\
X \ar[r]^i & Y
}
\]
In this situation, both $\Til X$ and $\Til Y$ are nonsingular.
Let $N$ be the normal bundle to $X$ in $Y$, of rank~$d$; identify 
$\Til X$ with the projectivization of $N$, and let $\zeta=c_1(\cO_N(1))$. 
Then (Theorem~15.4 in \cite{MR85k:14004}):
\[
c(\Til Y)-f^* c(Y) = j_*(g^* c(X)\cdot \alpha)\quad,
\]
where
\[
\alpha = \frac 1\zeta \left[
\sum_{i=0}^d g^* c_{d-i}(N) - (1-\zeta) \sum_{i=0}^d (1+\zeta)^i
g^* c_{d-i}(N)\right]\quad.
\]
Proofs of this formula that do not use Riemann-Roch were found by
A.~T.~Lascu and D.~B.~Scott (\cite{MR0360583}, \cite{MR0472816}).
In \cite{MR0472816}, Lascu and Scott write: {\em ``In this paper we
give a simple (and we hope definitive) proof of the result using only
simple arguments with vector bundles and some straightforward
manipulations.''\/}

\subsection{}
In this paper we go one step beyond the work of Lascu and Scott,
and prove `by simple arguments' a somewhat stronger result than the
formula recalled in \S\ref{introlead}. Our general aim is to remove the
nonsingularity hypothesis on $X$ and $Y$; this is what we accomplish,
pushing the level of generality to that of any {\em regular embedding\/}
of schemes $X\subseteq Y$.

As long as $X$ is regularly embedded in $Y$, the blow-up of $Y$ along 
$X$ can be regularly embedded into a projective 
bundle\footnote{As in \cite{MR85k:14004}, $P(E)$ denotes the projective 
bundle of lines in the vector bundle $E$.} $P(E)$ over $Y$ (see~\S\ref{ups}):
\[
\Til Y \subseteq P(E)\quad.
\]
The main result of this paper is the computation of the Chern classes
of the normal bundle $N_{\Til Y}$ of this embedding. In case~$X$ and
$Y$ (and hence $\Til X$ and $\Til Y$) are nonsingular, the Chern classes
of $T_{\Til Y}$ are immediately computed from the classes of $N_{\Til Y}$
and of $T_{P(E)}$; this recovers the formula recalled in \S\ref{introlead}
(cf.~\S\ref{oldrec}). The new proof of this formula appears to us simpler
than either the approach via Riemann-Roch or the proof found by 
Lascu and Scott.

\subsection{}\label{exer}
Let $Y$ be a scheme (pure dimensional, separated, of finite type over a field,
and admitting a closed embedding into a nonsingular scheme).

To state the result, assume first that $X$ is a complete intersection in $Y$: 
that is, $X$ is the zero-scheme of a regular section of a bundle $\hat N$ 
{\em on $Y$\/} of rank equal to the codimension of $X$. 
(Thus, the normal bundle $N$ to $X$ in $Y$ is isomorphic to the restriction 
of $\hat N$ to~$X$.) Let $f: \Til Y \to Y$ be the blow-up of $Y$ along $X$, 
and let $\Til X$ be the exceptional divisor.
Let $E$ be a vector bundle on $Y$, containing $\hat N$, 
and let $\hat C$ be the quotient:
\[
\xymatrix@1{
0 \ar[r] & \hat N \ar[r] & E \ar[r] & \hat C \ar[r] & 0
}\quad.
\]
There are natural embeddings $\Til Y\hookrightarrow P(\hat N)
\hookrightarrow P(E)$, which are easily seen to be regular 
(Lemma~\ref{reginci}); and $\cO(\Til X)$ is realized as the restriction of 
the universal subbundle $\cO(-1)$, hence $\cO(\Til X)\subseteq 
f^*(\hat N)$.  We prove:

\begin{lemma}\label{simlem}
With notation as above, let $N_{\Til Y}$ be the normal bundle of 
$\Til Y$ in $P(E)$. Then there is an exact sequence
\[
\xymatrix@1{
0 \ar[r] & f^*(\hat N)/\cO(\Til X) \ar[r] & N_{\Til Y} \ar[r] &
f^*(\hat C)\otimes \cO(-\Til X) \ar[r] & 0\quad.
}
\]
In particular,
\[
c(N_{\Til Y}) = \frac{c(f^*(\hat N))\, c(f^*(\hat C)\otimes \cO(-\Til X))}
{c(\cO(\Til X))}\quad.
\]
\end{lemma}

\subsection{}\label{gene}
In general, let $X\hookrightarrow Y$ be a regular embedding, but not 
necessarily a complete intersection. It is still the case 
(cf.~\cite{MR85k:14004}, B.8.2) that $X$ can be expressed as the 
zero-scheme of a section of a bundle $E$ on $Y$, and there is an 
exact sequence of vector bundles {\em on $X$:\/}
\[
\xymatrix{
0 \ar[r] & N \ar[r] & E|_X \ar[r] & C \ar[r] & 0
}
\]
where $N$ is the normal bundle of $X$ in $Y$, $E|_X$ denotes the 
restriction of $E$ to $X$, and $C$ is the quotient. The blow-up $\Til Y$
of $Y$ along $X$ still embeds regularly in $P(E)$ (\S\ref{ups}). 
The general result is as follows; parsing the formula requires some 
considerations, which follow the statement.

\begin{theorem}\label{main}
With notation as above, let $N_{\Til Y}$ be the normal bundle of $\Til Y$ in 
$P(E)$. Then
\[
c(N_{\Til Y})=\frac{c(\Nbb)\, c(\Cbb\otimes \cO(-\Til X))}
{c(\cO(\Til X))}\quad,
\]
where $\Nbb$, $\Cbb$ evaluate to the the pull-backs of $N$, resp.~$C$.
\end{theorem}

\subsection{Parsing}\label{parsing}
The terms $\Nbb$, $\Cbb$ appearing in the statement should be
understood as indeterminates with respect to which the right-hand-side
can be expanded; the terms in the expansion can then be interpreted
(by relating $\Nbb$ to $N$ and $\Cbb$ to $C$), determining a well-defined
operator on the Chow group $A_*(\Til Y)$. The content of the theorem
is that this operator equals $c(N_{\Til Y})$.

Here are the details of this operation. Expanding the expression gives
\[
\frac{c(\Nbb)\, c(\Cbb \otimes \cO(-\Til X))}{c(\cO(\Til X))}=
c(\Nbb)\, c(\Cbb) + \left( \cdots \right)\zeta
\]
where $\zeta=c_1(\cO(-\Til X))$, and the term $(\cdots)\, \zeta$ collects
monomials $c_i(\Nbb)\, c_j(\Cbb)\,\zeta^k$ with $k\ge 1$.
We prescribe that the first term should act on a class $a\in A_*(\Til Y)$
as the pull-back of $E$:
\[
c(\Nbb)\, c(\Cbb)\,(a) := f^*c(E)\cap a\quad.
\]
As for the remaining terms $c_i(\Nbb)\, c_j(\Cbb)\,\zeta^k$: if 
$a\in A_*\Til Y$, then $\zeta^k\cap a$ is supported on~$\Til X$ for $k\ge 1$, 
and hence $c_i(g^*N)\, c_j(g^*C) \zeta^k\cap (a)$ makes sense as a class 
in $A_*(\Til X)$, and determines (by proper push-forward) 
a class in $A_*(\Til Y)$. We prescribe
\[
c_i(\Nbb)\, c_j(\Cbb)\,\zeta^k (a) := 
j_*(c_i(g^*N)\, c_j(g^*C) \zeta^k\cap (a))\quad.
\]

In a nutshell, we want to think of $\Nbb$ and $\Cbb$ as pull-backs of
make-believe extensions to $Y$ of $N$ and $C$. If $N$ happens to
be the restriction of a bundle $\hat N$ (as in the complete intersection
case), $\hat N\subseteq E$, and $\hat C=E/\hat N$, then setting
$\Nbb=f^*(\hat N)$ and $\Cbb=f^*(\hat C)$ leads to the formula
presented above. It is a lucky circumstance that the formula can
be given a meaning even when $N$ is {\em not\/}
the restriction of a bundle defined on $Y$, and an even luckier
circumstance that the interpreted formula still computes the Chern
class of the normal bundle to the blow-up in the ambient projective
bundle.

\subsection{}\label{oldrec}
In the particular case when $X$ and $Y$ are nonsingular, 
Theorem~\ref{main} implies the formula recalled in~\S\ref{introlead}.
To see this, note that
\[
c(T_{P(E)}|_{\Til Y})=c(f^*E\otimes \cO(1))\, c(f^*T_Y)
\]
if $Y$ is nonsingular, by standard facts (for example, see B.5.8 in 
\cite{MR85k:14004}). Using the same parsing convention as in the
statement of Theorem~\ref{main}, this equals
\[
c(\Nbb \otimes \cO(1))\, c(\Cbb \otimes \cO(1))\, c(f^*T_Y)
\quad,
\]
and we get (with $\zeta=c_1(\cO(1))$):
\[
c(\Til Y)=\frac{c(T_{P(E)}|_{\Til Y})}{c(N_{\Til Y})}\cap [\Til Y]
=\frac{(1-\zeta)\, c(\Nbb \otimes \cO(1))}{c(\Nbb)}\cap f^* c(Y)\quad.
\]
This formula still uses the same convention: expand
\[
\frac{(1-\zeta)\, c(\Nbb \otimes \cO(1))}{c(\Nbb)}=
1+ \zeta(\cdots)\quad;
\]
replacing $c_i(\Nbb)$ by $g^*c_i(N)$ as explained above and capping
against $f^* c(Y)$ gives a class in $A_*(\Til Y)$. It is now easy to check
that this recovers on the nose the formula given in \S\ref{introlead}. 
The push-forward $j_*$ is responsible for the extra factor~$-1/\zeta$.

\subsection{}
If $Y$ is singular, but still a local complete intersection in a nonsingular
ambient variety $M$, then it admits a `virtual tangent bundle' 
$T^{\text{vir}}_Y$ (defined in $K$-theory as the difference between 
the restriction of $T_M$ and the normal bundle of $Y$ in $M$, see B.7.6 
in \cite{MR85k:14004}). Thus, $Y$ has well-defined Chern classes
$c(Y):=c(T^{\text{vir}}_Y)$. In this case $X$ and $\Til Y$ are also local 
complete intersections, and it is an easy consequence of Theorem~\ref{main}
that the formula given in \S\ref{introlead} holds if one uses these virtual
Chern classes throughout.

However, Theorem~\ref{main} is more general than this statement,
since it poses no restrictions on how singular $Y$ may be.

There are other notions of `Chern classes for singular varieties',
generalizing the nonsingular case, such as the 
Chern-Schwartz-MacPherson class $\csm$ mentioned below. 
It would be valuable to have formulas controlling the behavior of 
these classes under blow-ups at the level of generality considered
in this paper.

\subsection{}
As mentioned above, the proof of Theorem~\ref{main} appears to us
simpler than other approaches to the classical (and less general) formula.
Lemma~\ref{simlem} is a straightforward exercise; the extension 
from the complete intersection case to the general case follows from a 
mild generalization of a standard computational tool, namely B.6.10 
in \cite{MR85k:14004}.
On the other hand, it is worth noting that this generalization (proved 
in~\S\ref{B610}) ultimately relies on the technique known as 
{\em deformation to the normal cone;\/}
this is the technical tool behind the proofs found by Lascu and Scott,
as well as one of the main approaches to the proof of Riemann-Roch.
In fact, the reader may want to compare the `short' version of the proof
given in \S\ref{shortv}, which assumes familiarity with the deformation 
to the normal cone, and the detailed version given in \S\ref{longv}. 
The details in \S\ref{longv} are just as demanding as in the paper of 
Lascu and Scott.

Thus, it may be argued that these proofs of the blowing-up Chern
class formula are all different variations on
the same theme. Theorem~\ref{main} is a variation that happens to
work under the only requirement that $X$ be regularly embedded in~$Y$.

\subsection{}
The blowing-up Chern class formula has been used for calculations in
string theory (see for example~\cite{andreas-1999-3}); however, some
of my physicists acquaintances have expressed the opinion that the form 
recalled here in \S\ref{introlead} is difficult to apply, and its proof through
Riemann-Roch is somewhat obscure. I will close this introduction by
giving two short independent proofs of important particular cases,
which to my knowledge are not available in the literature. 
The formulation given in Lemma~\ref{difflp} may be more user-friendly
than the formula given in \S\ref{introlead}.

This subsection is independent of the rest of the paper, and (unlike the
rest) is limited to the case in which $X$ and $Y$ are nonsingular.

\subsubsection{Complete intersection, nonsingular case}\label{comintnsc}

\begin{lemma}\label{difflp}
Let $X\subseteq Y$ be nonsingular varieties.
If $X$ is a complete intersection of $d$ nonsingular hypersurfaces 
$Z_1,\dots, Z_d$ meeting transversally in $Y$, then
\begin{equation*}
\tag{*}
c(T_{\Til Y})=\frac{(1+\Til X)(1+f^*Z_1-\Til X)\cdots (1+f^* Z_d-\Til X)}
{(1+f^*Z_1)\cdots (1+f^* Z_d)}\cdot f^* c(T_Y)\quad.
\end{equation*}
\end{lemma}

\begin{proof}
By hypothesis, $Z=Z_1\cup\cdots \cup Z_d$ is a divisor with simple normal 
crossings in~$Y$, and it is easily checked that the divisor $W$ consisting of 
the exceptional divisor~$\Til X$ and of the proper transforms $W_i$ of 
$Z_i$ is a divisor with simple normal crossings in~$\Til Y$. We therefore 
have bundles
of tangent fields with logarithmic zeros (dual to the bundle of differential
forms with logarithmic poles) $T_Y(-\log Z)$, resp.~$T_{\Til Y}(-\log W)$ 
on $Y$, resp.~$\Til Y$. Comparing sections shows that
$T_{\Til Y}(-\log W) = f^* T_Y(-\log Z)$,
and hence
\[
c(T_{\Til Y}(-\log W)) =f^* c( T_Y(-\log Z))
\]
by the functoriality of Chern classes. Chern classes of bundles of tangent
fields with logarithmic zeros are well-known (see e.g.~Lemma~3.8 in
\cite{MR2183846}); we get
\[
\frac{c(T_{\Til Y})}{(1+\Til X)(1+f^*Z_1-\Til X)\cdots (1+f^* Z_d-\Til X)}
=\frac{f^* c(T_Y)}{(1+f^*Z_1)\cdots (1+f^* Z_d)}\quad,
\]
from which (*) follows immediately.
\end{proof}

Expanding formula~(*), and keeping in mind that $c_i(N_X)$ is the $i$-th
elementary symmetric function in $Z_1,\dots,Z_d$, one gets precisely the 
terms in the standard formulation presented in \S\ref{introlead}. In this sense, 
while Lemma~\ref{difflp} has a more limited scope ($X$ has to be a complete 
intersection), (*) may serve as mnemonics for the classical general formula,
and has a completely transparent proof.

\begin{remark}
Lemma~\ref{difflp} is a particular case of the following interesting fact.
Let $Z:=\sum Z_i$ be a divisor with normal crossings and nonsingular
components $Z_i$ in a nonsingular variety~$Y$. Say that a blow-up 
of $Y$ is `adapted to $Z$' if its center is the intersection of any
collection of the components $Z_i$. It is easily checked that in this
case the
exceptional divisor, together with the proper transforms of the
components of $Z$, form a divisor with simple normal crossings in
the blow-up. Say that a sequence of blow-ups over $Y$ is `adapted
to $Z$' if the first blow-up is adapted to $Z$, the second is adapted
to the new normal crossing divisor, etc.

Arguing as in Lemma~\ref{difflp}, one sees that if $\pi: \Til Y \to Y$
is any adapted sequence of blow-ups with respect to any divisor
$Z$ with simple normal crossings in $Y$, then $\pi^* \csm(\one_U)=
\csm(\one_{\Til U})$, where $U$ is the complement of $Z$ in $Y$
and $\Til U$ is the complement of $\pi^{-1}(Z)$ in $\Til Y$,
and $\csm$ denotes the `Chern class for constructible functions'
discussed below. (The $\csm$ class of the complement
of a normal crossing divisor is computed by the Chern class of the
tangent bundle with logarithmic zeros along the divisor, see 
e.g.~\cite{MR2001d:14008}, \S2.)
By standard functoriality properties of $\csm$, this formula holds as soon
as $\pi: \Til Y \to Y$ is a proper map dominated by a sequence of
adapted blow-ups. See \cite{AM}, \S4, for a more extensive
discussion, and for an application.
\end{remark}

\subsubsection{Characteristic zero, nonsingular case}
Over an algebraically closed field of characteristic zero, the formula of 
\S\ref{introlead} admits a very quick proof, without the complete 
intersection hypothesis of \S\ref{comintnsc}, if one takes for granted 
the theory of Chern classes for (possibly) singular varieties developed 
by Robert MacPherson\footnote{These classes are known to coincide, 
{\em mutatis mutandis,\/} with the classes defined earlier by 
M.-H.~Schwartz, see \cite{MR83h:32011}.} in \cite{MR50:13587}; see 
\cite{MR85k:14004}, \S19.1.7 for a version adapted to the Chow group.
According to this theory, 
there are `Chern classes' defined for every constructible function on a variety,
such that the Chern class of the constant function~$1$ on a nonsingular
variety equals the total Chern class of the tangent bundle. These
classes are covariant with respect to a push-forward of constructible
functions defined by taking Euler characteristics of fibers. The theory
is developed in characteristic~$0$; the basic covariance property does
not extend to positive characteristic (see~\S5.2 in~\cite{MR2282409}).

In the case of a blow-up map $f: \Til Y \to Y$ of a nonsingular variety
$Y$ along a codimension~$d$ nonsingular subvariety $X$, the Euler 
characteristic\footnote{This is the conventional topological Euler 
characteristic if the ground field is~$\Cbb$, and a suitable adaptation
over other algebraically closed fields of characteristic~zero.}
of the fibers is
\[
\chi(f^{-1}(p))=\left\{\aligned
1 &\quad p\not\in X \\
d &\quad p \in X
\endaligned
\right.\quad;
\]
it follows that 
\[
f_*(\one_{\Til Y}) = \one_Y + (d-1) \one_X\quad,
\]
where $\one$ denotes the constant function~1 on the given locus.
The covariance of Chern classes proved by MacPherson implies then
\begin{equation*}
\tag{1}
f_*(c(T_{\Til Y})\cap [\Til Y])=c(T_Y)\cap [Y] + (d-1)\, i_* c(T_X)\cap [X]
\quad,
\end{equation*}
where $i$ is the inclusion $X\hookrightarrow Y$.

On the other hand, it is easy to evaluate the restriction of $c(T_{\Til Y})$ to 
$\Til X$: 
\begin{equation*}
\tag{2}
j^* c(T_{\Til Y})=(1+\Til X)\, c(T_{\Til X})
=(1-\zeta)\, c(g^* N_X\otimes \cO(1))\, c(g^* T_X)\quad,
\end{equation*}
using the identification $\Til X\cong P(N_X)$, and with $\zeta=\cO(-\Til X)$. 

\begin{lemma}\label{check}
The class $c(T_{\Til Y})$ is characterized by formulas~(1) and~(2).
\end{lemma}

Indeed, {\em every\/} class in the Chow group of $\Til Y$ is characterized
by its push-forward to~$Y$ and its restriction to~$\Til X$ 
(\cite{MR85k:14004}, Proposition~6.7 (d)). It is now a simple exercise
(left to the reader)
to check that the formula for $c(\Til Y)$ stated~in \S\ref{introlead} 
satisfies both~(1) and~(2),
and by~Lemma~\ref{check} this suffices to prove the blowing-up
Chern class formula.

\subsection{}
In \cite{MR0412185}, Lascu and Scott propose a simplification of the
blow-up formula of~\S\ref{introlead}, that is equivalent to the formula
given here in Lemma~\ref{difflp}. However, they obtain this simpler formula
as a corollary of their blow-up formula; the proof of Lemma~\ref{difflp}
given in \S\ref{comintnsc} is independent (and essentially immediate).

In \cite{MR2363418},
Hansj\"org Geiges and Federica Pasquotto extend the classic blow-up
formula of \S\ref{introlead} to the case of symplectic,
complex, and real manifolds; their method follows closely the proof 
of Lascu and Scott in \cite{MR0472816}, whose algebro-geometric
ingredients they transfer to the topological environment.
\vskip 12pt

\subsection{}
{\em Acknowledgments.} I thank the Max-Planck-Institut in Bonn for 
hospitality and support. This work was also supported by NSA grant 
H98230-07-1-0024.


\section{Proof of Lemma~\ref{simlem}}

\subsection{}
We use notation as in \S\ref{exer}: $X$, $Y$ are pure dimensional
separated schemes of finite type over a field; $Y$ admits a closed
embedding into a nonsingular scheme. We assume that 
$X$ is a complete intersection in $Y$
of codimension~$d$, given as a the zero-scheme of a regular section 
of a vector bundle~$\hat N$ of rank~$d$; $f: \Til Y \to Y$ is the blow-up
of $Y$ along $X$. An exact sequence 
\[
\xymatrix@1{
0 \ar[r] & \hat N \ar[r] & E \ar[r] & \hat C \ar[r] & 0
}\quad.
\]
of vector bundles is given on $Y$. The embedding 
$\hat N\hookrightarrow E$ gives an embedding of projective bundles
\[
\xymatrix@C=5pt{
P(\hat N) \ar@{^(->}[rr] \ar[dr]_\pi & & P(\hat E) \ar[dl] \\
& Y
}
\]
with normal bundle $\pi^*(\hat C)\otimes \cO(1)$.

The section of $\hat N$ defining $X$ corresponds to a map
\[
\cO \hookrightarrow \cN
\]
to the sheaf of sections of $\hat N$; dualizing this map gives a surjection
\[
\cN^\vee \twoheadrightarrow \cI
\]
onto the ideal sheaf $\cI$ of $X$ in $Y$. Taking Proj of Sym gives
an embedding
\[
\Proj(\Sym^* \cI) \hookrightarrow P(\hat N)\quad.
\]
Now $\Sym^* \cI$ equals the Rees algebra of $\cI$, since the embedding
of $X$ in $Y$ is regular. Thus, $\Proj(\Sym^* \cI)=\Til Y$, and we have
fiberwise linear embeddings
\[
\xymatrix@1{
\Til Y \ar@{^(->}[r]^-\iota & P(\hat N) \ar@{^(->}[r] & P(E)\quad.
}
\]

\begin{lemma}\label{reginci}
The embedding $\iota$ is regular. 
\end{lemma}

\begin{proof}
The matter is local, so we may assume that $Y=\Spec A$, and that
the map $\cN^\vee \to \cO$ corresponds to a map $A^{\oplus d}
\to A$, where $d$ is the codimension of $X$ in $Y$. By assumption
$X$ is regularly embedded in $Y$, hence the image of this map is
an ideal $(a_1,\dots,a_d)$ generated by a regular sequence in $A$.
The blow-up $\Til Y$ is defined  by the equations $a_i T_j-a_j T_i$,
$1\le i<j\le d$, in $Y\times \Pbb^{d-1}=Y\times P(\hat N)$ 
(\cite{MR85k:14004}, Lemma~A.6.1). On the open set $U_d$ of 
$P(\hat N)$ defined by $T_d\ne 0$, the ideal of the blow-up is
\[
(a_1-a_d x_1,\dots,a_{d-1}-a_d x_{d-1})\quad,
\]
where $x_i=T_i/T_d$. Thus $\Til Y\cap U_d$ is a complete intersection 
in $U_d$. The situation is of course analogous on all open charts
$U_k=\{T_k\ne 0\}$. The statement follows.
\end{proof}

Note that the tautological line bundle $\cO(-1)$ on $P(E)$ restricts
to its namesakes on $P(\hat N)$, on $\Til Y=\Proj(\Sym^* \cI)$, and
on the exceptional divisor $\Til X=P(N_XY)$. Further, 
$\cO(\Til X)\cong \cO(-1)|_{\Til Y}$; this determines an embedding
of $\cO(\Til X)$ in $\iota^*\pi^* (\hat N)=f^*(\hat N)$.

\subsection{}
At this point we have maps as in the commutative diagram:
\[
\xymatrix{
\Til Y \ar@{^(->}[r]^-\iota \ar[rd]_f & P(\hat N) \ar@{^(->}[r] \ar[d]^\pi
& P(E) \ar[ld] \\
& Y
}
\]
The regular embeddings $\Til Y \hookrightarrow P(\hat N) \hookrightarrow
P(E)$ yield an exact sequence of normal bundles
\[
\xymatrix{
0 \ar[r] & N_{\Til Y} P(\hat N) \ar[r] & 
N_{\Til Y} P(E) \ar[r] & 
\iota^* N_{P(\hat N)} P(E) \ar[r] & 0\quad.
}
\]
Letting $N_{\Til Y}$ denote $N_{\Til Y}P(E)$ as in Lemma~\ref{simlem},
this is
\[
\xymatrix{
0 \ar[r] & N_{\Til Y} P(\hat N) \ar[r] & 
N_{\Til Y} \ar[r] & 
f^*(\hat C) \otimes \cO(1) \ar[r] & 0
}\quad,
\]
and in order to prove Lemma~\ref{simlem} it suffices to verify the
following:

\begin{lemma}\label{preplem}
\[
N_{\Til Y} P(\hat N) \cong f^*(\hat N)/\cO(\Til X)\quad.
\]
\end{lemma}

\subsection{Proof of Lemma~\ref{preplem}.}
Let $\cK$ be the kernel of the surjection $\cN^\vee \to \cI$.
Taking $\Sym$, we obtain the exact sequence
\[
\xymatrix{
0 \ar[r] & \cK \cdot \Sym^{*-1} \cN^\vee \ar[r] &
\Sym^* \cN^\vee \ar[r] & 
\oplus_{k\ge 0} \cI^k \ar[r] & 0
}
\]
determining the ideal of $\Til Y$ in $P(\hat N)$; it follows that the
conormal sheaf to $\Til Y$ in $P(\hat N)$ is $f^* \cK \otimes \cO(-1)$.

Pulling back to $\Til Y$ the sequence $0 \to \cK \to \cN^\vee \to
\cO_Y \to \cO_X \to 0$, we get the sequence
\begin{equation*}
\tag{*}
\xymatrix{
0\ar[r] & { f^*(\cK)} \ar[r] & f^*(\cN^\vee) \ar[r] & \cO_{\Til Y} \ar[r] & 
\cO_{\Til X} \ar[r] & 0\quad,
}
\end{equation*}
and I claim that replacing $f^*(\cK)$ by $f^*(\cK)\otimes \cO(-1)$ in this
sequence produces an {\em exact\/} sequence on $\Til Y$:
\begin{equation*}
\tag{**}
\xymatrix{
0\ar[r] & { f^*(\cK)}\otimes \cO(-1) \ar[r] & f^*(\cN^\vee) \ar[r] & \cO_{\Til Y} 
\ar[r] & \cO_{\Til X} \ar[r] & 0\quad.
}
\end{equation*}
Indeed, use again the local presentation obtained in 
the proof of Lemma~\ref{reginci}: we start from the exact sequence
\[
\xymatrix{
0 \ar[r] & (a_i T_j-a_j T_i)_{1\le i<j\le d} \ar[r] & 
A^{\oplus d} \ar[r] & A \ar[r] & A/(a_1,\dots,a_d) \ar[r] & 0
}\]
where $T_i$ is the generator of the $i$-th factor in the middle, and
$A^{\oplus d} \to A$ is defined by $T_i \mapsto a_i$; the fact that
the kernel is as stated follows from the fact that $(a_1,\dots,a_d)$
is regular. Pull-back to a representative chart in the blow-up by 
tensoring by
\[
B=\frac{A[x_1,\dots,x_{d-1}]}{(a_1-a_d x_1,\dots,a_{d-1}-
a_d x_{d-1})}\quad:
\]
this yields the sequence corresponding to (*):
\[
\xymatrix{
0 \ar[r] & (a_d(T_i-x_i T_d))_{1\le i< d} \ar[r] & 
B^{\oplus d} \ar[r] & B \ar[r] & B/(a_d) \ar[r] & 0
}
\]
The morphism $B^{\oplus d} \to B$ is still defined by $T_i \mapsto a_i$,
and its kernel is easily checked to be $(T_i-x_i T_d)_{1\le i< d}$
($a_d$ is a non-zero-divisor in $B$). We see 
that this differs from $f^* \cK$ by the presence in the latter of the 
extra factor of $a_d$. As $a_d$ is a section of $\cO(1)$, dividing by 
$a_d$ corresponds to tensoring by $\cO(-1)$, and this concludes the 
verification that the sequence (**) is exact.

Now rewrite (**) as the exact sequence of locally free sheaves on $\Til Y$:
\[
\xymatrix{
0\ar[r] & f^*(\cK)\otimes \cO(-1) \ar[r] & f^*(\cN^\vee) \ar[r] & \cO(-\Til X)
\ar[r] & 0\quad.
}
\]
Dualizing, and using the identification of $f^*(\cK)\otimes \cO(-1)$ with
the conormal sheaf to $\Til Y$ in $P(\hat N)$, gives the exact sequence of 
vector bundles on $\Til Y$:
\[
\xymatrix{
0\ar[r] & \cO(\Til X) \ar[r] & f^*(\hat N) \ar[r] & N_{\Til Y} P(\hat N)
\ar[r] & 0\quad,
}
\]
concluding the proof of Lemma~\ref{preplem}.\hfill\qed

As noted above, Lemma~\ref{simlem} follows from Lemma~\ref{preplem}.
Thus, we have now established Theorem~\ref{main} under the hypothesis 
that $X$ is a complete intersection in $Y$.


\section{Proof of Theorem~\ref{main}}\label{proofmain}

\subsection{}\label{ups}
Let $X\hookrightarrow Y$ be a regular embedding. As recalled in 
\S\ref{gene}, we can express
$X$ as the zero-scheme of a section of a bundle $\rho: E\to Y$, and we 
have an exact sequence of vector bundles on~$X$:
\[
\xymatrix{
0 \ar[r] & N \ar[r] & E|_X \ar[r] & C \ar[r] & 0
}
\]
where $N$ is the normal bundle of $X$ in $Y$. The blow-up $\Til Y$
of $Y$ along $X$ embeds in~$P(E)$, and this embedding is regular: 
indeed, this is a local matter, so it reduces to the case considered 
in Lemma~\ref{reginci}.

We view this situation as follows. The section of $E$ defining $X$ is
an embedding
\[
s: Y \hookrightarrow E
\]
of $Y$ into the total space of $E$. Both $s$ and the zero-section 
$z: Y \hookrightarrow E$ are regular embeddings, with normal
bundle $E$ itself; and $s(Y)$, $z(Y)$ meet along $X$. In other
words, we have the fiber square
\begin{equation*}
\tag{$\dagger$}
\xymatrix@1{
X \ar@{^(->}[r] \ar@{^(->}[d] & s(Y) \ar@{^(->}[d] \\
z(Y) \ar@{^(->}[r] & E
}
\end{equation*}
in which all embeddings are regular.

In particular, this gives an embedding of the normal 
bundle $N$ to $X$ in $Y=s(Y)$ into the restriction of the
normal bundle to $z(Y)$ in $E$, that is, $E|_X$.

\subsection{}
Now we let $\nu: \Til E\to E$ be the blow-up along $z(Y)$. 
The blow-up $\Til Y=B\ell_XY$ may be realized as the proper
transform of $s(Y)$ in $\Til E$. 

Note that $z(Y)$ is a complete intersection in $E$: it is the zero-scheme 
of the `identity' section $E \to \rho^*(E)$. Thus, we are in the situation
of Lemma~\ref{simlem}, and we can conclude that $\Til E$ embeds 
regularly into $P(\rho^*(E))$, with normal bundle
\[
N_{\Til E}P(\rho^*(E))\cong \frac{\nu^*(E)}{\cO(\Til W)}\quad,
\]
where $\Til W$ denotes the exceptional divisor. Summarizing, we have
the commutative diagram
\[
\xymatrix@1{
\Til Y \ar@{^(->}[r] \ar@{^(->}[d] & P(E) \ar@{^(->}[d] \\
\Til E \ar@{^(->}[r] & P(\rho^*(E))
}
\]
of regular embeddings, where the vertical map on the left is the proper
transform of~$s(Y)$ (this will be verified to be a regular embedding in
Lemma~\ref{rempt}), and the vertical map on the right is obtained
by restricting $P(\rho^*(E))$ to $s(Y)$. It follows that
\[
c(N_{P(E)} P(\rho^*(E)))\cdot c(N_{\Til Y} P(E))
=c(N_{\Til E} P(\rho^*(E)))\cdot c(N_{\Til Y}\Til E)\quad,
\]
(omitting pull-backs for convenience), and hence
\[
c(N_{\Til Y})=\frac{c(N_{\Til E} P(\rho^*(E)))\cdot c(N_{\Til Y}\Til E)}
{c(N_{P(E)} P(\rho^*(E))}
\]
where $N_{\Til Y}=N_{\Til Y}P(E)$ as in Theorem~\ref{main}.
Since $N_{\Til E}P(\rho^*(E))\cong \nu^*(E)/\cO(\Til W)$ by 
Lemma~\ref{simlem}, and $N_{P(E)}P(\rho^*(E))$ is
the pull-back of $N_{s(Y)}E$, that is $E$, and further
$\cO(\Til W)$ restricts to $\cO(\Til X)$ on $\Til Y$, 
we can conclude that
\[
c(N_{\Til Y})=\frac{c(N_{\Til Y}\Til E)}{c(\cO(\Til X))}\quad,
\]
reducing the computation of $c(N_{\Til Y})=c(N_{\Til Y}P(E))$, which 
is our objective, to the computation of $c(N_{\Til Y}\Til E)$.

\subsection{}\label{concluproof}
The proof of Theorem~\ref{main} is now complete if we show:

\begin{claim}\label{finn}
With the notational convention explained in \S\ref{parsing},
\[
c(N_{\Til Y}\Til E)=c(\Nbb)\, c(\Cbb\otimes \cO(-\Til X))\quad.
\]
\end{claim}
This is an instance of a general template, which appears to be
independently useful, and which we treat in the next section. 
Claim~\ref{finn} is the result of applying Theorem~\ref{newnormal} 
to the situation of diagram~($\dagger$). Therefore, the proof of 
Theorem~\ref{newnormal} will conclude the proof of 
Theorem~\ref{main} (and this paper).


\section{The normal bundle of a proper transform}\label{B610}

\subsection{}
Let $X\subseteq Y$ and $Y\subseteq Z$ be regular embeddings, and let
$\Til Y$, $\Til Z$ be the blow-ups along $X$; $\Til Y$ may be identified
with the proper transform of $Y$ in $\Til Z$. Then
(\cite{MR85k:14004}, B.6.10) $\Til Y$ is regularly embedded in $\Til Z$,
and
\[
N_{\Til Y}\Til Z\cong f^*(N_YZ)\otimes \cO(-\Til X)\quad,
\]
where $f:\Til Y \to Y$ is the blow-up map, and $\Til X$ is the exceptional 
divisor in $\Til Y$. 

We wish to extend this formula (at the level of Chern classes) to the case
in which the center $W$ of the blow-up is not necessarily contained in $Y$,
but $X=W\cap Y$ is still regularly embedded in both $W$ and $Y$: all
embeddings in the diagram
\begin{equation*}
\tag{$\ddagger$}
\xymatrix@1{
X \ar@{^(->}[r] \ar@{^(->}[d] & Y \ar@{^(->}[d] \\
W \ar@{^(->}[r] & Z
}
\end{equation*}
are regular. Note that the diagram $(\dagger)$ of \S\ref{proofmain}
is an instance of this situation: take $Y=s(Y)$, $W=z(Y)$, $Z=E$.

\subsection{}\label{TYTZ}
Let $\Til Z \to Z$ be the blow-up along~$W$. The blow-up $\Til Y$
of $Y$ along $W\cap Y=X$ embeds in $\Til Z$ as the proper transform
of $Y$.

\begin{lemma}\label{rempt}
$\Til Y$ is regularly embedded in $\Til Z$.
\end{lemma}

\begin{proof}
This is a local verification, which follows closely the case $X=W$ proved
in \cite{MR85k:14004}, B.6.10 (from which we already borrowed in the
proof of Lemma~\ref{reginci}). We may assume that $Z=\Spec A$, 
the ideal of $W\subset Z$ is generated by a regular sequence
$(a_1,\dots,a_d)$, and the ideal of $Y\subset Z$ is also generated
by a regular sequence $(a_1,\dots,a_e,b_1,\dots, b_\ell)$, with 
$1\le e\le d$.
The blow-up $\Til Z$ is defined by $a_i T_j-a_j T_i$,
$1\le i<j\le d$, in $Z\times \Pbb^{d-1}=Z\times P(\hat N)$. On the open 
set defined by $T_d\ne 0$, $\Til Z$ has coordinate ring
\[
\Til A:=\frac{A[x_1,\dots,x_{d-1}]}{(a_1-a_d x_1,\dots,a_{d-1}-a_d x_{d-1})}\quad,
\]
where $x_i=T_i/T_d$.

At the same time, $Y$ has coordinate ring $A'=A/(a_1,\dots,a_e,b_1,\dots, b_\ell)$;
by assumption, the cosets $\overline a_{e+1},\dots,\overline a_d\in A'$ of
$a_1,\dots,a_d$ form a regular sequence. The coordinate ring of a matching
chart for $\Til Y$ is{\small
\begin{multline*}
\Til A':=\frac{A'[x_{e+1},\dots,x_{d-1}]}{(\overline a_{e+1}-\overline a_d x_{e+1},\dots,
\overline a_{d-1}-\overline a_d x_{d-1})} \\
\cong \frac{A[x_{e+1},\dots,x_{d-1}]}{(b_1,\dots,b_\ell,a_1,\dots,a_e,
a_{e+1}-a_d x_{e+1},\dots, a_{d-1}- a_d x_{d-1})}
\quad.
\end{multline*}}
On this chart, the inclusion $\Til Y\subset \Til Z$ corresponds to the surjection
$\Til A \twoheadrightarrow \Til A'$ of $A$-algebras given by $x_1\mapsto 0,
\dots,x_e\mapsto 0,\dots, x_{e+1}\mapsto x_{e+1},\dots, x_d\mapsto x_d$.
The kernel of this surjection is generated by
\[
x_1,\dots,x_e,b_1,\dots,b_\ell\quad,
\]
clearly a regular sequence at each point of $\Til Y$.
 This verifies that the embedding is regular on
this chart, and the situation is identical in the other charts $T_k\ne 0$,
$k>e$. (The argument also implies that $\Til Y$ has empty intersection
with the charts $T_k\ne 0$, $k\le e$.)
\end{proof}

The challenge is to compute $c(N_{\Til Y}\Til Z)$. At one extreme, 
$X=W$ and we are in the situation of \cite{MR85k:14004}, B.6.10:
in this case $c(N_{\Til Y}\Til Z)=c(f^*(N_YZ)\otimes \cO(\Til X))$.

At the other extreme, $Y$ and $W$ intersect properly in $Z$:
that is, $\Til Y$ equals the total transform of $Y$ in $\Til Z$; 
in this case, $c(N_{\Til Y}\Til Z)=c(f^*(N_YZ))$.

The general case lies `in between' these two special cases.

\subsection{}\label{partc}
In the fiber square $(\ddagger)$:
\[
\xymatrix@1{
X \ar@{^(->}[r] \ar@{^(->}[d]_i & Y \ar@{^(->}[d] \\
W \ar@{^(->}[r] & Z
}
\]
note that there is an embedding $N_XY\subset i^*N_WZ$, 
and therefore an exact sequence
\[
\xymatrix{
0 \ar[r] & N_XY \ar[r] & i^* N_WZ \ar[r] & C \ar[r] & 0
}
\]
The cokernel $C$ is the {\em excess normal bundle\/} of the square
(cf.~\cite{MR85k:14004}, \S6.3).

As a useful warm-up, assume that there is a regularly embedded 
subscheme $Z'$ of~$Z$ containing $Y$ and $W$, and in which 
$Y$ and~$W$ intersect properly:
\[
\xymatrix@1@C=10pt@R=5pt{
X \ar@{^(->}[rr] \ar@{^(->}[dd]_i & & Y \ar@{^(->}[dd] \\
\\
W \ar@{^(->}[rr] & & Z' \ar@{^(->}[dr] \\
& & & Z
}
\]
Also, assume that all embeddings are regular; and let $h:\Til Z'\to Z'$ be the
blow-up of $Z'$ along $W$. 

Denote by $\Nbb$ the normal bundle of $Y$ in $Z'$, and its pull-backs;
and denote by $\Cbb$ the normal bundle of $Z'$ in $Z$, as well as its 
pull-backs. Note that
\[
c(N_YZ)=c(\Nbb)\, c(\Cbb)\quad,
\]
while the fact that $Y$ and $W$ meet properly in $Z'$ implies that
$i^* N_WZ'=N_XY$, and hence that $\Cbb$ restricts to $C$ on $X$.
By the same token, $\Nbb$ restricts to $N_XW$ on $X$.

This situation is a combination of the two `extremes' mentioned at the 
end of \S\ref{TYTZ}:

---Since $Y$ and $W$ meet properly in $Z'$, we have 
\[
c(N_{\Til Y}\Til Z')=c(f^*N_YZ')=c(\Nbb)\quad;
\]

---Since $Z'$ contains the center $W$ of the blow-up, we have
\[
c(N_{\Til Z'}\Til Z)=c(h^*N_{Z'}Z\otimes \cO(1))=c(\Cbb\otimes \cO(1))
\quad,
\]
where $\cO(-1)$ stands for the line bundle of the exceptional divisor. 

---Therefore,
\[
c(N_{\Til Y}\Til Z)=c(\Nbb)\, c(\Cbb\otimes \cO(1))
\]
(omitting evident pull-backs).

Our main result is that this formula holds in the general case
(even if $Z'$ is not present), provided that it is interpreted appropriately.

\subsection{The statement}\label{thestat}
Summarizing: in general, two bundles are defined on $X$, namely
$N_XW$ and the excess intersection bundle $C$. In the particular
case considered in \S\ref{partc}, these two bundles extend to 
bundles $\Nbb$, resp.~$\Cbb$ defined on the whole of $Y$, such
that $c(N_YZ)=c(\Nbb) c(\Cbb)$, and we have verified that
\[
c(N_{\Til Y}\Til Z)=c(\Nbb)\, c(\Cbb\otimes \cO(\Til X))
\]
where pull-backs via $f: \Til Y \to Y$ are understood.

Here is how the right-hand-side of this formula may be interpreted
as an operator on $A_*\Til Y$, even when $N_XW$ and $C$ are
not assumed to be restrictions of bundles $\Nbb$, $\Cbb$ 
(cf.~\S\ref{parsing}).
\[
\xymatrix{
\Til X \ar[r]^j \ar[d]_g & \Til Y \ar[d]^f \\
X \ar[r]^i & Y
}
\]
\begin{itemize}
\item[---] Formally expand $c(\Nbb)\, c(\Cbb\otimes \cO(\Til X))$:
\[
c(\Nbb)\, c(\Cbb\otimes \cO(\Til X))=c(\Nbb)\, c(\Cbb) + \PP(c_i(\Nbb),
c_j(\Cbb))\cdot \Til X
\]
for a well-defined polynomial $\PP$ in the (formal) variables $c_1(\Nbb),
c_2(\Nbb),\dots$, $c_1(\Cbb), c_2(\Cbb), \dots$;
\item[---] For $\alpha\in A_*\Til Y$, define
\[
\boxed{
c(\Nbb)\, c(\Cbb\otimes \cO(\Til X))\cap \alpha:=
f^* c(N_YZ)\cap \alpha + j_* \PP(c_i(g^*N_XW), c_j(g^*C))\cap (\Til X
\cdot \alpha)}
\quad.
\]
\end{itemize}

\begin{theorem}\label{newnormal}
With notation as above,
\[
c(N_{\Til Y}{\Til Z})\cap \alpha=c(\Nbb)\, c(\Cbb\otimes \cO(\Til X))
\cap \alpha
\]
for all $\alpha\in A_*\Til Y$.
\end{theorem}

In the application to $(\dagger)$ in \S\ref{proofmain}, $N_XW$ equals
the normal bundle $N$ of $X$ in (the image via the zero-section of) $Y$,
and $C$ equals the cokernel of the inclusion of $N$ into $E|_X$, as 
prescribed in \S\ref{gene}. Thus, Theorem~\ref{newnormal} does provide
the last ingredient in the proof of Theorem~\ref{main}, as pointed out
in \S\ref{concluproof}. Proving Theorem~\ref{newnormal} is our last task.

\subsection{Proof of Theorem~\ref{newnormal}: short version}\label{shortv}
The following summary will suffice for the expert. The deformation to the
normal cone (\cite{MR85k:14004}, Chapter~5) may be used to reduce the 
general situation $(\ddagger)$ to the `linearized' situation
\begin{equation*}
\tag{$\ddagger'$}
\xymatrix@1{
X \ar@{^(->}[r] \ar@{^(->}[d] & N_XY \ar@{^(->}[d] \\
N_XW \ar@{^(->}[r] & N_XZ
}
\end{equation*}
This is covered by the particular case considered in \S\ref{partc}, by taking
$Z'$ to be the (direct) sum of $N_XY$ and $N_XW$ in $N_XZ$. As 
shown in \S\ref{partc} the formula holds in this case, hence it holds in
general.

We end this article by spelling out this argument, for the benefit of readers
who may be less familiar with the deformation to the normal cone.

\subsection{Proof of Theorem~\ref{newnormal}: long version}\label{longv}
The following diagram may be helpful in tracing the argument:
\[
\xymatrix@C=5pt{
\Til Y\cong \Til Y\times \{0\} \ar@{^(->}[rr]^-z && \Til M_Y
:=B\ell_{X\times \Pbb^1} M_Y \ar@/^1pc/[rrrd]^-p \ar[dr]^-\nu \ar[dl]_\varphi \\
& M_Y :=B\ell_{X\times \{\infty\}} Y\times \Pbb^1 \ar[dr] && 
\Til Y\times \Pbb^1 \ar[dl] \ar[rr] && \Til Y\\
&& Y\times \Pbb^1 \ar[d] \\
&& \Pbb^1
}
\]
Here $M_Y$ is the deformation of $Y$ to the normal cone (bundle) $N_XY$.
The subscheme $X\times \Pbb^1$ of $Y\times \Pbb^1$ lifts to an isomorphic
copy in $M_Y$, and $\varphi: \Til M_Y\to M_Y$ is the blow-up along
this isomorphic copy. It is easily checked that the inverse image of 
$X\times \Pbb^1\subset Y\times \Pbb^1$ in $\Til M_Y$ is a Cartier divisor, 
and more precisely
it equals the sum of the two exceptional divisors; by the universal 
property of blow-ups, the map $\Til M_Y \to \Til Y\times \Pbb^1$ factors
through $\Til Y\times\Pbb^1$, as indicated in the diagram. In fact, 
$\nu: \Til M_Y \to \Til Y\times \Pbb^1$ is the blow-up along $\Til X\times
\{\infty\}$.

Also, note that the composition
\[
p\circ z\quad:\quad \Til Y \to \Til M_Y \to \Til Y
\]
is the identity. 

With $M_W$, resp.~$M_Z$ obtained similarly from $W\times \Pbb^1$, 
resp.~$Z\times \Pbb^1$ by blowing up along $X\times \{\infty\}$, 
we have inclusions
\[
\xymatrix@1{
X\times \Pbb^1 \ar@{^(->}[r] \ar@{^(->}[d] & M_Y \ar@{^(->}[d] \\
M_W \ar@{^(->}[r] & M_Z
}
\]
Over all $t\ne \infty$, this diagram specializes to $(\ddagger)$; 
over $\infty$, the diagram formed by the exceptional divisors:
\begin{equation*}
\tag{$\ddagger''$}
\xymatrix@1{
X \ar@{^(->}[r] \ar@{^(->}[d] & P(N_XY\oplus 1) 
\ar@{^(->}[d] \\
P(N_XW\oplus 1) \ar@{^(->}[r] & P(N_XZ\oplus 1)
}
\end{equation*}
is the projective completion of the `linearized' diagram $(\ddagger')$.
At $\infty$ we also find copies of $\Til Y$, $\Til W$, $\Til Z$, meeting
the corresponding projective completions along their exceptional
divisors. The scheme-theoretic intersection of $M_Y$ and 
$M_W$ is the lift of $X\times \Pbb^1$ mentioned above. This locus
is disjoint from the copy of $\Til Y$ at $\{\infty\}$.
Further, $\varphi$ restricts to an isomorphism of the proper transform
via $\nu$ of $\Til Y\times\{\infty\}$ (which is isomorphic to $\Til Y$
as $\Til X$ is a divisor in $\Til Y$) with this copy of $\Til Y$ at $\infty$
in $M_Y$.

Blow-up $M_Z$ along $M_W$; the proper transform of $M_Y$
agrees with the blow-up of the latter along $X\times \Pbb^1$, so it
is the variety $\Til M_Y$ appearing in the larger diagram.
Over any $t\neq \infty$ (and in particular for $t=0$), the blow-ups
reproduce the situation considered in~\S\ref{TYTZ}. 

We have to verify that $c(N_{\Til Y}\Til Z)\cap \alpha=c(\Nbb) 
c(\Cbb\otimes \cO(\Til X))\cap \alpha$ for all $\alpha\in A_*\Til Y$. 
Letting $\Gamma\cdot \alpha:=j_* \PP(c_i(g^*N_XW), c_j(g^*C))\cap 
\Til X\cdot \alpha$ as in the definition preceding the statement of 
Theorem~\ref{newnormal}, the task is to show that
\[
(c(N_{\Til Y}\Til Z)-c(f^*N_YZ))\cap \alpha = \Gamma\cdot \alpha
\]
for all $\alpha\in A_*\Til Y$, and we have verified that this holds in the
situation considered in~\S\ref{partc}. We let $\Gamma^*$ be the 
operator defined in the same way as $\Gamma$ on $A_*(\Til M_Y$), 
and observe that $\Gamma^*$ restricts to $\Gamma$ over all $t\neq \infty$, 
and to the analogous operator for the linearized version ($\ddagger''$).

By linearity, we may assume that $\alpha=[V]$, where $V\subset \Til Y$
is a subvariety of $\Til Y$.
Since $N_{\Til Y}\Til Z$, resp.~$f^* N_YZ$ may be realized as
pull-backs via~$z$ of $N_{\Til M_Y}\Til M_Z$, resp.~$\varphi^*
N_{M_Y}M_Z$, the projection formula gives
\begin{equation*}
\tag{*}
\small{
(c(N_{\Til Y}\Til Z)-c(f^*N_YZ))\cap [V] = 
p_*\left((c(N_{\Til M_Y}\Til M_Z)-c(\varphi^* N_{M_Y}M_Z))
\cap ([V\times \{0\}])\right)\quad.}
\end{equation*}

The proper transform of $V\times \Pbb^1\subset \Til Y\times 
\Pbb^1$ in $\Til M_Y$ is the blow-up $M_V$ along
$(\Til X\cap V)\times\{\infty\}$; the fiber of $M_V$ over
$\{0\}$ is precisely the variety $V\times \{0\}$ appearing in~(*).
This is rationally equivalent to the fiber of $M_V$ over
$\{\infty\}$, that is
\[
P(N_{\Til X\cap V}V\oplus 1)\cup B\ell_{\Til X\cap V}V\quad.
\]
Thus, 
\begin{multline*}
(c(N_{\Til Y}\Til Z)-c(f^*N_YZ))\cap [V] \\ = 
p_*\left((c(N_{\Til M_Y}\Til M_Z)-c(\varphi^* N_{M_Y}M_Z))
\cap ([P(N_{\Til X\cap V}V\oplus 1)]+[B\ell_{\Til X\cap V}V])\right)
\quad.
\end{multline*}
As noted earlier, $\varphi$ restricts to an isomorphism from
$B\ell_{\Til X\cap V}V\cong V$ to $V\subset \Til Y\subset M_Y$.
The target $V$ is disjoint from the center $X\times \Pbb^1$ of the
blow-up~$\varphi$, therefore
\[
(c(N_{\Til M_Y}\Til M_Z)-c(\varphi^* N_{M_Y}M_Z))
\cap [B\ell_{\Til X\cap V}V]=0 \quad,
\]
and hence
\[
(c(N_{\Til Y}\Til Z)-c(f^*N_YZ))\cap [V] \\ = 
p_*\left((c(N_{\Til M_Y}\Til M_Z)-c(\varphi^* N_{M_Y}M_Z))
\cap [P(N_{\Til X\cap V}V\oplus 1)]\right)\,.
\]
Now we are squarely in the blow-up over the linearized
diagram $(\ddagger'')$. This situation is contemplated by the
case considered in \S\ref{partc}: use $P(N_XY\oplus N_XW
\oplus 1)$ for $Z'$. Therefore, the theorem holds in this case,
giving
\[
(c(N_{\Til Y}\Til Z)-c(f^*N_YZ))\cap [V] \\ = 
p_*\left(\Gamma^* \cdot [P(N_{\Til X\cap V}V\oplus 1)]\right)
\quad.
\]
Next we essentially run through the construction in reverse.
Since $\Gamma^*$ is supported on the exceptional divisor
of $\varphi$, $\Gamma^*\cdot [B\ell_{\Til X\cap V}V]=0$, hence
\[
(c(N_{\Til Y}\Til Z)-c(f^*N_YZ))\cap [V] \\ = 
p_*\left(\Gamma^* \cap( [P(N_{\Til X\cap V}V\oplus 1)]+
[B\ell_{\Til X\cap V}V])\right)
\quad;
\]
since $[P(N_{\Til X\cap V}V\oplus 1)]+[B\ell_{\Til X\cap V}V]
=[V\times\{0\}]$ in $M_V$, 
\[
(c(N_{\Til Y}\Til Z)-c(f^*N_YZ))\cap [V] \\ = 
p_*\left(\Gamma^* \cdot [V\times \{0\}]\right)
\quad;
\]
and since $\Gamma^*$ restricts to $\Gamma$ on fibers over
$t\ne \infty$, the projection formula gives
\[
(c(N_{\Til Y}\Til Z)-c(f^*N_YZ))\cap [V] = \Gamma\cdot [V]
\]
as claimed. This concludes the proof of Theorem~\ref{newnormal}.
\hfill\qed



\end{document}